\documentstyle[bezier,12pt]{article}

\input{amssym.def}	
\input{amssym}		
\newcommand{\mathbb}[1]{{\Bbb #1}}
\newcommand{\mthrm}[1]{\mbox{#1}}
\newcommand{\ensuremath}[1]{#1}

\newcommand{\qed}{$\;\;\;\Box$}
\newenvironment{proof}{\par\smallbreak{\sl Proof.~}}
{\unskip\nobreak\hfill \qed \par\medbreak}
\newcommand{\hide}[1]{}

\newcommand{\Con}{\ensuremath{\mthrm{C}}}
\newcommand{\D}{\ensuremath{{\cal D}}}
\newcommand{\E}{\ensuremath{{\cal E}}}

\newcommand{\mb}[1]{\ensuremath{\mathbb{#1}}}
\newcommand{\N}{\mb{N}}
\newcommand{\R}{\mb{R}}

\newcommand{\A}{\ensuremath{{\cal A}}}
\newcommand{\G}{\ensuremath{{\cal G}}}

\newcommand{\EM}{\ensuremath{{\cal E}_{\mthrm{M}}}}
\newcommand{\NN}{\ensuremath{{\cal N}}}

\renewcommand{\d}{\ensuremath{\partial}}
\newcommand{\grad}{\ensuremath{\mbox{\rm grad}\,}}

\newtheorem{thm}{Theorem}

\newtheorem{defn}[thm]{Definition}

\newtheorem{rem}[thm]{Remark}
\newtheorem{ex}[thm]{Example}

\newcommand{\al}{\alpha}

\newcommand{\ga}{\gamma}

\newcommand{\eps}{\varepsilon}

\newcommand{\om}{\omega}

\newcommand{\supp}{\mathop{\mthrm{supp}}}
\newcommand{\diag}{\mathop{\rm diag}\nolimits}

\title{
Generalized Solutions to Hyperbolic Systems with Nonlinear
Conditions and Strongly Singular Data
}

\newcounter{thesame}
\author{
I.~Kmit\\
{\small
Institute for Applied Problems of Mechanics and Mathematics,}\\
{\small Ukrainian Academy of Sciences}\\
{\small Naukova St.\ 3b}\\
{\small 79060 Lviv,
Ukraine}
\\
{\small   E-mail:
{\tt kmit@ov.litech.net}}
}

\date{}

\begin{document}

\maketitle

\begin{abstract}
Using  the framework of Colombeau algebras of generalized
functions, we prove the existence and uniqueness results for global
generalized solvability of semilinear
hyperbolic systems with nonlinear
 nonlocal boundary conditions.
We admit
strong singularities in the differential equations as well as in
the initial and boundary conditions.
Our analysis covers the case of
non-Lipshitz nonlinearities both in the differential equations
and in the boundary conditions.

\end{abstract}

\section{Introduction}\label{sec:intr}

We study existence and uniqueness  of global
generalized solutions to mixed problems for semilinear hyperbolic
systems with nonlinear nonlocal boundary conditions.
Specifically,
in the domain $\Pi=\{(x,t)\,|\,0<x<l$, $t>0\}$
we study the following problem:
\begin{eqnarray}
(\partial_t  + \Lambda(x,t)\partial_x)
 U &=& F(x,t,U), \qquad (x,t)\in \Pi
  \label{eq:1}\\[3mm]
U(x,0) &=& A(x), \qquad x\in (0,l) \label{eq:2} \\[3mm]
U_i(0,t) &=& H_i(t,V(t)),\quad k+1\le i\le n, \qquad  t\in(0,\infty)\nonumber\\
U_i(l,t) &=& H_i(t,V(t)),\quad 1\le i\le k, \qquad  t\in(0,\infty) \label{eq:3} \, ,
\end{eqnarray}
where  $U$, $F$, and $A$  are real $n$-vectors,
$\Lambda=\diag(\Lambda_1,\dots,
\Lambda_n)$ is a diagonal matrix,
$\Lambda_1,\dots,\Lambda_k<0$, $\Lambda_{k+1},\dots,\Lambda_n>0$
for some $1\le k\le n$, and
$V(t)=(U_1(0,t),\dots,$
$U_{k}(0,t),U_{k+1}(l,t),\dots,U_{n}(l,t))$.
Due to the conditions imposed on
$\Lambda$, the system~(\ref{eq:1}) is non-strictly
hyperbolic. Note also that the boundary of $\Pi$ is not characteristic.
We will denote $H=(H_1,\dots,H_n)$.

Special cases of~(\ref{eq:1})--(\ref{eq:3}) arise in
laser dynamics \cite{r1,rw,sbw,tlo} and
chemical kinetics~\cite{z}.

All the data of the problem are allowed to be
strongly singular, namely, they can be of any desired order of
singularity. This entails nonlinear superpositions for
(strongly singular)
distributions in the right-hand sides of ~(\ref{eq:1})--(\ref{eq:3})
(including
compositions of the singular initial data and the singular
characteristic curves). To tackle this complication, we use
the framework of Colombeau
algebra of generalized functions $\G(\overline\Pi)$~\cite{bia,col,Obe}.
We show that all superpositions appearing here
are well defined in $\G(\overline\Pi)$.

We establish a positive existence-uniqueness result in $\G(\overline\Pi)$
for the problem~(\ref{eq:1})--(\ref{eq:3}) with strongly singular initial
data and with nonlinearities of the following type (more detailed description
is given in Section~3):
The functions $F$ and $H$ are either
Lipshitz with Colombeau generalized numbers as Lipshitz constants or
non-Lipshitz with less than   quadratic growth
in $U$ and $V$.

For different aspects of the subject we refer the reader to
sources~\cite{Gr1,Gr2,hoop,HKm,LOb,Obe,Obe2,NP}. The essential
assumption made on $F$ in papers~\cite{HKm,Obe} is that $\grad_UF$ is
globally bounded uniformly over $(x,t)$ varying in any compact set.
The main complication with the
non-Lipshitz nonlinearities when investigating Colombeau solutions
lies in the following.
The
solutions in $\G$ are nets of smooth functions which are classical solutions
to the associated problems with smooth initial data. To guarantee the
existence of such solutions, in general, one has to assume that the
nonlinearities have bounded gradients.
In the present paper we tackle this complication combining classical and
nonclassical approaches.

Papers~\cite{NP,NP1} deal with Cauchy
problems for semilinear hyperbolic systems~(\ref{eq:1}) with
$F$ slowly increasing at infinity.
The nonlinear
term is replaced by a suitable regularization $F_{\eps}$
having a bounded gradient
with respect to $U$ for every fixed $\eps$ and converging to
$F$ as $\eps\to 0$.
The regularized system is solved in $\G(\R^2)$.
Moreover, in~\cite{NP} the
components of $\Lambda$ are allowed to be 1-tempered generalized
functions. The authors replace $\Lambda$ by its regularization
(the regularization procedure is similar to the one for $F$) which is
1-tempered generalized function of bounded growth and solve the regularized
problem.

 In~\cite{Gr1,Gr2} the author investigates
weak limits for
 semilinear hyperbolic systems and nonlinear superpositions
for strongly singular distributions appearing in these systems.
  He establishes an optimal
link between the singularity of the initial data and the growth of
the nonlinear term.
Weak limits of strongly singular Cauchy problems for semilinear
hyperbolic systems with bounded, sublinear, and superlinear
growth are investigated
in~\cite{dem-ra,Km1,mo1,ra-re}.

Existence-uniqueness results within Colombeau algebras for two-dimensional
hyperbolic problems with discontinuous coefficients $\Lambda_i$
are obtained in~\cite{hoop,LOb,Obe,Obe2}.
These papers impose an essential restriction on the coefficients,
which allows one to avoide
the negative effect of infinite propagation speed.
Namely, the coefficients are assumed to be
globally bounded in the Colombeau algebra $\G$.
At the present paper
we do not assume global boundedness of
the $\Lambda_i$ in~(\ref{eq:1}), thereby allowing them to be strongly
singular (this is also the case for~\cite{Alu}).

A novelty of the paper is that it treats (within $\G(\overline\Pi)$)
strongly singular initial data
of the problem (including coefficients $\Lambda_i$
in~(\ref{eq:1})),
nonlinear boundary conditions,
and non-Lipshitz nonlinearities in~(\ref{eq:1}) and~(\ref{eq:3}).

The plan of our exposition is as follows. In Section~2 we compile some
facts about Colombeau algebra of generalized functions.
In Section~3 we state and prove our main results.
We prove the existence-uniqueness result within
$\G(\overline\Pi)$ for Colombeau Lipshitz nonlinearities
(Subsection 3.1) and extend its to the case of
non-Lipshitz nonlinearities (Subsection 3.2).

\section{Preliminaries}

In this section we summarize the relevant
material on the full version of Colombeau algebras of generalized functions.

Let $\Omega\subset\R^n$ be a domain in $\R^n$. We denote by $\G(\Omega)$
and $\G(\overline\Omega)$ the full version of Colombeau algebra of
generalized functions over $\Omega$ and $\overline\Omega$, respectively.
To define $\G(\Omega)$ and $\G(\overline\Omega)$, we
first introduce the mollifier spaces used to parametrize
the regularizing sequences of generalized functions.
Given $q\in \N_0$ denote
$$
\begin{array}{cc}
\displaystyle
\A_q(\R)=\Bigl\{\varphi\in\D(\R)\Bigm|\int\varphi(x)\,dx=1,\int x^k\varphi(x)\,dx=0
\mbox{\ for \ } 1\le k\le q\Bigr\},\\[7mm]
\A_q(\R^n)=\Bigl\{\varphi(x_1,\dots,x_n)=\prod\limits_{i=1}^n\varphi_0(x_i)\Bigm|
\varphi_0\in\A_q(\R)\Bigr\}.
\end{array}
$$
If $\varphi\in\A_0(\R^n)$, let
$$
\varphi_{\varepsilon}(x)=\frac{1}{\varepsilon^n}\varphi\biggl(
\frac{x}{\varepsilon}\biggr).
$$
Set
$$
\E(\overline{\Omega})=\{u: \A_0\times\overline{\Omega}\to \R \Bigm|
u(\varphi,.)\in\Con^{\infty}(\overline{\Omega})\ \  \forall \varphi\in\A_0(\R)\}.
$$
We define the algebra of moderate elements
$\EM(\overline{\Omega})$ to be  the
subalgebra of $\E(\overline{\Omega})$  consisting of the
elements $u\in\E(\overline{\Omega})$
such that
$$
\begin{array}{cc}
\forall K\subset \overline{\Omega} \mthrm{\ compact}, \forall \alpha\in\N_0^n,
\exists N\in\N \mthrm{\ such\  that \ } \forall \varphi\in\A_N(\R^n)
 \\[5mm]
\displaystyle
\exists C>0, \exists \eta>0
\mthrm{\ with \ }
\sup\limits_{x\in K}|\d^{\alpha}u(\varphi_{\varepsilon},x)|\le C
\varepsilon^{-N}, \ \ 0<\varepsilon<\eta.
\end{array}
$$
The ideal $\NN(\overline{\Omega})$ consists of all $u\in\EM(\overline{\Omega})$
such that
$$
\begin{array}{cc}
\forall K\subset \overline{\Omega} \mthrm{\ compact}, \forall \alpha\in\N_0^n,
\exists N\in\N \mthrm{\ such\  that \ } \forall q\ge N,
\forall \varphi\in\A_q(\R^n)
 \\[5mm]
\displaystyle
\exists C>0, \exists \eta>0
\mthrm{\ with \ }
\sup\limits_{x\in K}|\d^{\alpha}u(\varphi_{\varepsilon},x)|\le C\varepsilon^{q-N}, \ \ 0<\varepsilon<\eta.
\end{array}
$$
Finally,
$$
\G(\overline{\Omega})=\EM(\overline{\Omega})/\NN(\overline{\Omega}).
$$
This is an associative and commutative differential algebra.
The algebra $\G({\Omega})$ on open set is constructed in the same manner,
with $\Omega$ in place of $\overline{\Omega}$. Note that
$\G({\Omega})$ admits a
canonical embedding of $\D'(\Omega)$. We will use the notation
$U=[(u(\varphi,x))_{\varphi\in\A_0(\R^n)}]$ for the elements $U$ of
$\G(\Omega)$ with representative $u(\varphi,x)$.

One of the advantages of using Colombeau algebra of generalized functions $\G$
lies in the fact that in a variety of important cases
the division by generalized functions, in particular the
division by discontinuous functions and measures, is defined in $\G$.
Complete description of the cases when the division is possible
in the full version of Colombeau algebras
is given by the following criterion of
invertibility~\cite{Alu} (the criterion of invertibility for special
version of Colombeau algebras $\G_s(\Omega)$ is proved in~\cite{book}):

\begin{thm}\label{thm:invert}
Let $U\in\G(\Omega)$ ($U\in\G(\overline\Omega)$). Then the following two
conditions are equivalent:\\
(i) $U$ is invertible in $\G(\Omega)$
(in $\G(\overline\Omega)$), i.e., there exists
$V\in\G(\Omega)$ ($V\in\G(\overline\Omega)$) such that $UV=1$ in $\G(\Omega)$
(in $\G(\overline\Omega)$).\\
(ii) For each representative $(u(\varphi,x))_{\varphi\in\A_0(\R^n)}$ of $U$
and each compact set $K\subset\Omega$ ($K\subset\overline\Omega$)
there exists $p\in\N$ such that for all $\varphi\in\A_p(\R^n)$ there is
$\eta>0$ with $\inf\limits_{K}|u(\varphi_{\varepsilon},x)|
\ge\varepsilon^p$ for all $0<\varepsilon<\eta$.
\end{thm}

\section{Existence-uniqueness results in the Colombeau algebra of generalized
functions}

\subsection{Colombeau Lipshitz nonlinearities}

We here develop some results of~\cite{Alu} and~\cite{HKm}  to the case
of nonlinear nonlocal boundary conditions and Colombeau Lipshitz nonlinearities
in~(\ref{eq:1}) and~(\ref{eq:3}) (with Lipshitz constants
as Colombeau generalized numbers).
We will need a notion of a
generalized function whose growth is more restrictive
than the
$1/\varepsilon$-growth (as in the definition of $\EM$).

\begin{defn}\label{defn:ga}(\cite{Alu})
Let $\Omega\subset\R^n$ be a domain in
$\R^n$.
Suppose we have a function $\gamma : (0,1)\mapsto (0,\infty)$.
We say that an element $U\in\G(\Omega)$ ($U\in\G(\overline{\Omega})$) is
{\rm locally of $\gamma$-growth}, if it has a representative
$u\in\EM(\Omega)$ ($u\in\EM(\overline\Omega)$) with the following property:

For every compact set $K\subset\Omega$ ($K\subset\overline\Omega$)
there is $N\in\N$ such that for every $\varphi\in \A_{N}(\R^n)$
there exist $C>0$ and $\eta>0$ with
$\sup\limits_{x\in K}|u(\varphi_{\varepsilon},x)|\le C\gamma^{N}(\varepsilon)$ for
$0<\varepsilon<\eta$.

\end{defn}

We now make assumptions on the initial data
of the problem~(\ref{eq:1})--(\ref{eq:3}). Let $\gamma(\varepsilon)$ and
$\gamma_1(\varepsilon)$ be functions from
$(0,1)$ to $(0,\infty)$ such that
\begin{equation}\label{eq:ga}
{\gamma(\varepsilon)}^{\gamma^N(\varepsilon)}=
O\biggl(\frac{1}{\varepsilon}\biggr),
\quad {\gamma(\varepsilon)}^{\gamma_1^N
(\varepsilon)}
=O\biggl(\frac{1}{\varepsilon}\biggr)\quad \mbox{as}\,\,\varepsilon\to 0
\end{equation}
for each $N\in\N$. Assume that

\begin{enumerate}

\item
$\Lambda(x,t)\in(\G(\overline\Pi))^n$, $A(x)\in(\G[0,l])^n$.

\item
$\Lambda_i$  for $i\le n$
 are locally of $\gamma$-growth on $\overline\Pi$ and
invertible on $\overline\Pi$.

\item $\partial_x\Lambda_i$ for $i\le n$
 are locally of $\gamma_1$-growth on $\overline\Pi$.

\item
$F(x,t,y)\in(\G(\overline\Pi\times\R^n))^n$,
$H(t,z)\in(\G([0,\infty)\times\R^n))^n$.

\item
For every compact set $K\subset\overline{\Pi}$ and $i\le n$
the mapping $y\mapsto F_i(x,t,y)$ and all its derivatives
are polynomially bounded for all $(x,t)\in K$  with coefficients in $\G(K)$.

\item
For every compact set $K\subset[0,\infty)$ and $i\le n$
the mapping $z\mapsto H_i(t,z)$
  and all its derivatives are
polynomially   bounded for all $t\in K$
 with coefficients in $\G(K)$.

\item
For every compact set $K\subset\overline{\Pi}$ there exists a nonnegative
generalized function $L_F(x,t)\in\G(K)$ such that for all
$(x,t)\in K$, $i\le n$, and
$y^1,y^2\in\R^n$ we have
$$
|F_i(x,t,y^1)-F_i(x,t,y^2)|
\le L_F(x,t)\sum\limits_{j=1}^n |y_j^1-y_j^2|.
$$

\item
For every compact set $K\subset[0,\infty)$ there exists a nonnegative
generalized function $L_H(x,t)\in\G(K)$ such that for all
$t\in K$, $i\le n$, and $z^1,z^2\in\R^n$ we have
$$
|H_i(t,z^1)-H_i(t,z^2)|\le L_H(t)\sum\limits_{j=1}^n |z_j^1-z_j^2|.
$$

\item
$L_F(x,t)$ and $L_H(t)$
 are locally of $\gamma$-growth on $\overline\Pi$ and $[0,\infty)$,
respectively.

\item $\supp A_i(x)\subset (0,l)$ and $\supp H_i(t,0)
\subset(0,\infty)$ for $i\le n$.
\end{enumerate}

Assumptions imposed on $\Lambda_i$ allow them to be strongly singular
and, even more, to have any desired order of singularity.
Assumptions 4--6 state that, given
$U\in(\G(\overline\Pi))^n$ and
$V\in(\G[0,\infty))^n$, $F(x,t,U)$ and $H(t,V)$
are well defined in
Colombeau algebra $\G$. We can interpret Assumptions 7 and 8 as the
Lipshitz conditions in Colombeau sense imposed on generalized functions
$F$ and $H$. Assumption~9 allows $L_F$ and $L_H$ to be strongly
singular.
The last assumption ensures  the compatibility of~(\ref{eq:2})
and~(\ref{eq:3}) of any desired order.

We now  state the main result of this subsection.
\begin{thm}\label{thm:gen}
Under Assumptions 1--10 where the functions $\gamma$ and $\gamma_1$
are specified by~(\ref{eq:ga}),  the
problem~(\ref{eq:1})--(\ref{eq:3})
has a unique solution $U\in\G(\overline{\Pi})$.
\end{thm}

\begin{proof}
We will first prove the existence  of a classical smooth solution to
the problem~(\ref{eq:1})--(\ref{eq:3}) where  the initial data are smooth,
satisfy Assumption 10, and the functions $F$ and $H$ have bounded
gradients with respect to $U$ and $V$, respectively,
uniformly over $(x,t)$ varying in compact subsets of
$\overline{\Pi}$.
In parallel, we will obtain a priori estimates for
classical smooth solutions and their derivatives of any desired order.
Therewith we will obtain the existence of a prospective
representative $u$ of the solution $U$ in
$\G(\overline\Pi)$.  To finish
the existence part of the proof,  we will show
the moderateness of $u$. The uniqueness of the constructed
generalized solution will be proved by the same scheme.

We first reduce
the problem~(\ref{eq:1})--(\ref{eq:3}) with smooth initial data
 to an equivalent integral-operator form. Denote by
$\omega_i(\tau;x,t)$ the $i$-th characteristic of~(\ref{eq:1})
passing through a point $(x,t)\in\overline{\Pi}$. From the
assumptions imposed on $\Lambda$ it follows that  such characteristic exists,
is smooth in $\tau,x,t$, and can be continued up to the boundary of
$\Pi$. The smallest value of
$\tau\ge 0$ at which the characteristic
$\xi=\omega_i(\tau;x,t)$ intersects $\partial\Pi$
will be denoted by $t_i(x,t)$.
Integrating each equation of~(\ref{eq:1}) along the corresponding
characteristic curve,
we obtain the following equivalent integral-operator
form of~(\ref{eq:1})--(\ref{eq:3}):
\begin{equation}\label{eq:integral}
\begin{array}{cc}
\displaystyle
U_i(x,t)=(R_iU)(x,t)+
\int\limits_{t_i(x,t)}^t\Bigl[U(\omega_i(\tau;x,t),\tau)\int\limits_0^1\nabla_U
F_i(\omega_i(\tau;x,t),\tau,\sigma U)\,d\sigma
\\[5mm]\displaystyle
+F_i(\omega_i(\tau;x,t),\tau,0)\Bigr]\,d\tau,\ \
1\le i\le n,
\end{array}
\end{equation}
where
$$
(R_iU)(x,t)=
\cases{V_i(t_i(x,t))\int\limits_0^1\nabla_V
H_i(t_i(x,t),\sigma V)\,d\sigma+H_i(t_i(x,t),0)&if
$t_i(x,t)>0$,
 \cr
A_{i}(\omega_i(0;x,t)) &if
$t_i(x,t)=0$. \cr}
$$
Given $T>0$, denote
$$
\Pi^T=\{(x,t)\,|\,0<x<l, 0<t<T\}.
$$
Set
$$
E_U(\al_1,\al_2;T)=\max\Bigl\{|\partial_x^{\al_1}\partial_t^{\al_2} U_i(x,t)|
\,\Big|\,
(x,t)\in\overline{\Pi}^T,
1\le i\le n\Bigr\},
$$
$$
E_F(\al_1,\al_2)=\max\Bigl\{|
\partial_x^{\al_1}\partial_t^{\al_2}F_i(x,t,y)|\,\Big|\,
(x,t,y)\in\overline{\Pi}^T\times
\{y:|y|\le E_U(0,0;T)\},1\le i\le n\Bigr\},
$$
$$
E_H(\al)=\max\Bigl\{|
\partial_t^{\al}H_i(t,z)|\,\Big|\,
(t,z)\in[0,T]\times
\{z:|z|\le E_U(0,0;T)\},1\le i\le n\Bigr\},
$$
$$
L_F^{\max}=\max\Bigl\{L_F(x,t)\,\Big|\,(x,t)\in\overline\Pi^T\Bigr\},\quad
L_H^{\max}=\max\Bigl\{L_H(t)\,\Big|\,t\in[0,T]\Bigr\}.
$$
Simplifying the notation, we drop the dependence of
$E_F(\al_1,\al_2)$, $E_H(\al)$,
$L_F^{\max}$, and $L_H^{\max}$ on $T$. Note that we will
use these parameters for a fixed $T>0$.

Assume that the initial data $\Lambda$, $F$, $A$, and $H$ of our
problem  are smooth with respect to all their arguments,
satisfy Assumption~10, and the functions $F$ and $H$ have bounded
gradients with respect to $U$ and $V$, uniformly over $(x,t)$
varying in compact subsets of $\overline{\Pi}$.
Fix an arbitrary $T>0$.
If $(x,t)\in\overline\Pi^T$, then $L_F^{\max}$ and $L_H^{\max}$
are Lipshitz constants of $F$ and $H$ with respect to $U$ and
$V$, respectively.
We now prove that the problem~(\ref{eq:integral}) has  a smooth
solution in $\overline\Pi^T$.
In parallel, we obtain  global a priori estimates
for smooth solutions, we will make use of for construction
of a Colombeau solution. We obtain the global a priori estimates
by iterating the a apriori estimates for local smooth solutions
in a number of steps. The proof is split in
four claims.

{\it Claim 1. The problem~(\ref{eq:integral}) has  a unique  continuous
solution in $\overline\Pi^T$.}
We start from the local continuous solution to~(\ref{eq:integral}),
namely, we state that there exists a unique solution
$U\in(\Con(\overline\Pi^{t_0}))^n$ to the problem~(\ref{eq:integral})
for some   $t_0>0$.
To prove this, choose $t_0$ satisfying the condition
\begin{equation}\label{eq:22_0}
\om_n(t;0,\tau)<\om_1(t;l,\tau)\quad
 \forall \tau\ge 0, \forall t\in[\tau,\tau+t_m]
\end{equation}
with $m=0$.
For $t\in[0,t_0]$ we can express
$V(t)$ in the form
\begin{equation}\label{eq:C}
\begin{array}{cc}
\displaystyle
V_i(t)=A_{i}(\omega_i(0;0,t))+
\int\limits_{0}^t\Bigl[U(\omega_i(\tau;0,t),\tau)\int\limits_0^1\nabla_U
F_i(\omega_i(\tau;0,t),\tau,\sigma U)\,d\sigma
\\[5mm]\displaystyle
+F_i(\omega_i(\tau;0,t),\tau,0)\Bigr]\,d\tau,\ \
1\le i\le k,
\\[5mm]\displaystyle
V_i(t)=A_{i}(\omega_i(0;l,t))+
\int\limits_{0}^t\Bigl[
U(\omega_i(\tau;l,t),\tau)\int\limits_0^1\nabla_U
F_i(\omega_i(\tau;l,t),\tau,\sigma U)\,d\sigma
\\[5mm]\displaystyle
+F_i(\omega_i(\tau;l,t),\tau,0)\Bigr]\,d\tau,\ \
k+1\le i\le n.
\end{array}
\end{equation}
Since~(\ref{eq:integral})
 is a system of Volterra integral equations of the
second kind in $\overline\Pi^{t_0}$, we can apply the contraction mapping principle.
We apply the operator defined by the right hand side of (\ref{eq:integral})
to two
continuous functions $U^1$ and $U^2$.
 Note that these functions
have the same initial and boundary values.
We cosider their difference in
$\overline\Pi^{t_0}$.
Notice the estimate
$$
E_{U^1-U^2}(0,0;t_0)\le t_0q_0E_{U^1-U^2}(0,0;t_0),
$$
where
$$
q_0=nL_F^{\max}(1+nL_H^{\max}).
$$
We are able to choose $t_0$ so that the additional condition
$t_0<1/q_0$ is obeyed. Then the contraction property of the operator
defined by the right hand side of (\ref{eq:integral}) holds
with respect to $\overline\Pi^{t_0}$.
We have thus proved
existence and uniqueness of a continuous solution $U$ to the
problem~(\ref{eq:integral}) in $\overline\Pi^{t_0}$.
Furthermore, we have the
following local a priori estimate:
\begin{equation}\label{eq:24}
\begin{array}{ccccc}
\displaystyle
E_U(0,0;t_0)\le
\frac{1}{1-q_0t_0}
\biggl[\Bigl(\max\limits_{x\in[0,l],1\le i\le n}|A_{i}(x)|
\\
\displaystyle+
T\max\limits_{(x,t)\in \overline{\Pi}^T,1\le i\le n}|F_i(x,t,0)|
\Bigr)\Bigl(1+nL_H^{\max}\Bigr)+\max\limits_{t\in[0,T],1\le i\le n}
|H_i(t,0)|\biggr].
\end{array}
\end{equation}
Note that the value of $q_0$ depends on $T$ and does not depend on $t_0$.
This allows us to complete the proof of the claim
in $\lceil T/t_0\rceil$ steps,
iterating local existence-uniqueness result in domains
$$
(\Pi^{jt_0}\cap\Pi^T)\setminus\overline\Pi^{(j-1)t_0}, \quad
1\le j\le \lceil T/t_0\rceil.
$$
Moreover,
using the estimate~(\ref{eq:24}) $\lceil T/t_0\rceil$ times and
each time starting  with the final
value of $U$ from the previous
step, we derive the following global a priori estimate:
\begin{equation}\label{eq:240}
\begin{array}{ccccc}
\displaystyle
E_U(0,0;T)\le
P_{1,0}\biggl(\frac{1}{1-q_0t_0},n,L_H^{\max}\biggr)\\
\displaystyle
\times
P_{2,0}\biggl(\max\limits_{x\in[0,l],1\le i\le n}|A_{i}(x)|,
\max\limits_{(x,t)\in \overline{\Pi}^T,1\le i\le n}|F_i(x,t,0)|,
\max\limits_{t\in[0,T],1\le i\le n}|H_i(t,0)|\biggr),
\end{array}
\end{equation}
where $P_{1,0}$ is a polynomial of degree $3\lceil T/t_0\rceil$
with all coefficients identically equal to 1 and $P_{2,0}$
is a polynomial of the first degree with positive constant
coefficients depending only on $T$.

{\it Claim 2. The problem~(\ref{eq:1})--(\ref{eq:3}) has  a unique
$\Con^1$-solution in $\overline\Pi^T$.}
The proof is similar to the proof of Claim~1.
 Let us consider the initial-boundary problem for
$\d_xU$:
\begin{equation}\label{eq:U_x}
\begin{array}{cc}
\displaystyle
\d_xU_i(x,t)=(R_{ix}^{'}U)(x,t)
+\int\limits_{t_i(x,t)}^t\Bigl[\nabla_UF_i(\xi,\tau,U)\cdot\d_xU
\\\displaystyle-
(\d_x\Lambda_i)(\xi,\tau)\d_xU_i+(\d_xF_i)(\xi,\tau,U)\Bigr]\Big|_{\xi=\om_i(\tau;x,t)}
\,d\tau,\quad
1\le i\le n,
\end{array}
\end{equation}
where
$$
(R_{ix}^{'}U)(x,t)
=\cases{\Lambda_i^{-1}(0,\tau)\Bigl[
F_i(0,\tau,U)-\cr
\nabla_VH_i(\tau,V)\cdot V^{'}(\tau)-
(\d_tH_i)(\tau,V)\Bigr]\Big|_{\tau=t_i(x,t)}
&if
$t_i(x,t)>0,k+1\le i\le n$,
 \cr
\Lambda_i^{-1}(l,\tau)\Bigl[
F_i(l,\tau,U)-\cr
\nabla_VH_i(\tau,V)\cdot V^{'}(\tau)-
(\d_tH_i)(\tau,V)\Bigr]\Big|_{\tau=t_i(x,t)}
&if
$t_i(x,t)>0,1\le i\le k$,
 \cr
A_{i}^{'}(\om_i(0;x,t)) &if
$t_i(x,t)=0$. \cr}
$$
Choose $t_1$ satisfying the condition~(\ref{eq:22_0}) with $m=1$.
Combining~(\ref{eq:1})
with~(\ref{eq:U_x}),
 we get
\begin{equation}
\begin{array}{cc}
V_i^{'}(t)=F_i(0,t,U)-\Lambda_i(0,t)(\d_xU_i)(0,t)\\
\displaystyle
=F_i(0,t,U)-\Lambda_i(0,t)\biggl[A_{i}^{'}(\om_i(0;0,t))+
\int\limits_{0}^t\Bigl[\nabla_UF_i(\xi,\tau,U)\cdot\d_xU\\
\displaystyle
-(\d_x\Lambda_i)(\xi,\tau)\d_xU_i+(\d_xF_i)(\xi,\tau,U)\Bigr]
\Big|_{\xi=\om_i(\tau;0,t)}
\,d\tau\biggr],\quad 1\le i\le k,
\end{array}
\end{equation}
where $t\in[0,t_1]$.
The functions $V_i^{'}(t)$ for $k+1\le i\le n$ can be expressed
in the same form.
Using the fact that $U$ is a known continuous function (see Claim~1),
we now
 apply the operator defined by the right hand side of (\ref{eq:U_x})
to two continuous functions $\d_xU^1$ and $\d_xU^2$.
 Note that these functions
have the same initial and boundary values.
We cosider their difference in
$\overline\Pi^{t_1}$.
Notice the estimate
$$
E_{U^1-U^2}(1,0;t_1)\le t_1q_1E_{U^1-U^2}(1,0;t_1),
$$
where
$$
q_1=(nL_F^{\max}+E_{\Lambda}(1,0;T))(1+nL_H^{\max}).
$$
We are able to choose $t_1$ so that the additional condition
$t_1<1/q_1$ is obeyed.
This shows that the operator
defined by the right hand side of (\ref{eq:U_x}) has the
contraction property
with respect to the domain
$\overline\Pi^{t_1}$.
Thus, we have proved the existence and the uniqueness of a
solution $U\in\Con_{x,t}^{1,0}(\overline\Pi^{t(1)})$ to the
problem~(\ref{eq:integral}). Furthermore,
we have the
following local a priori estimate:
\begin{equation}\label{eq:U1}
\begin{array}{ccccc}
\displaystyle
E_U(1,0;t_1)\le
\frac{1}{1-q_1t_1}
\biggl[\Bigl(\max\limits_{x\in[0,l],1\le i\le n}|A_{i}^{'}(x)|
+TE_F(1,0)\\[5mm]
+E_{\Lambda^{-1}}(0,0;T)E_F(0,0)
\Bigr)\Bigl(1+nL_H^{\max}\Bigr)
+E_{\Lambda^{-1}}(0,0;T)E_H(1)\biggr].
\end{array}
\end{equation}
Using the fact that the value of $q_1$ depends on $T$ and does
not depend on $t_1$
and iterating the local existence-uniqueness result in domains
$$
(\Pi^{jt_1}\cap\Pi^T)\setminus\overline{\Pi}^{(j-1)t_1}, \quad
1\le j\le \lceil T/t_1\rceil,
$$
we obtain the global a priori estimate:
\begin{equation}\label{eq:U1g}
\begin{array}{ccccc}
\displaystyle
E_U(1,0;T)\le P_{1,1}\biggl(\frac{1}{1-q_1t_1},n,L_H^{\max}\biggr)\\
\displaystyle
\times
P_{2,1}\biggl(\max\limits_{x\in[0,l],1\le i\le n}|A_{i}^{'}(x)|,
\max\limits_{0\le\al\le 1}E_F(\al,0),E_{\Lambda^{-1}}(0,0;T),E_H(1)\biggr),
\end{array}
\end{equation}
where $P_{1,1}$ is a polynomial of degree $3\lceil T/t_1\rceil$
with all coefficients identically equal to 1 and $P_{2,1}$
is a polynomial of the second degree
with positive constant
coefficients depending only on $T$.

The a priori estimate for  $E_U(0,1;T)$ now follows from the system~(\ref{eq:1}):
$$
E_U(0,1;T)\le E_F(0,0)+
E_{\Lambda}(0,0;T)
E_U(1,0;T),
$$
where $E_U(1,0;T)$ satisfies the estimate~(\ref{eq:U1g}).
This finishes the proof of the claim.

{\it Claim 3. The problem~(\ref{eq:1})--(\ref{eq:3}) has  a unique
$\Con^2$-solution in $\overline\Pi^T$.}
Following the proof of Claims~1 and~2, let us
consider the following problem for
$\d_x^2U$:
\begin{equation}\label{eq:U_xx}
\begin{array}{cc}
\displaystyle
\d_x^2U_i(x,t)=(R_{ixx}^{''}U)(x,t)
+\int\limits_{t_i(x,t)}^t\Bigl[\nabla_UF_i(\xi,\tau,U)\cdot\d_x^2U
\\[7mm]-
2(\d_x\Lambda_i)(\xi,\tau)\d_x^2U_i
-(\d_x^2\Lambda_i)(\xi,\tau)\d_xU_i
+(\d_x^2F_i)(\xi,\tau,U)+2\nabla_U(\d_xF_i)(\xi,\tau,U)\cdot\d_xU
\\[5mm]
+
\nabla_U(\nabla_UF_i(\xi,\tau,U)\cdot\d_xU)\cdot\d_xU
\Bigr]\Big|_{\xi=\om_i(\tau;x,t)}
\,d\tau,
1\le i\le n,
\end{array}
\end{equation}
where
$$
(R_{ixx}^{''}U)(x,t)
=A_{i}^{''}(\om_i(0;x,t))\quad \mbox{if}\,\,
t_i(x,t)=0
$$
and
$$
(R_{ixx}^{''}U)(x,t)=
-(\d_t\Lambda_i^{-1})(0,\tau)(\d_xU_i)(0,\tau)-\Lambda_i^{-2}(0,\tau)
\biggl[(\d_tF_i)(0,\tau,U)
$$
$$
+\nabla_UF_i(0,\tau,U)\cdot(\d_tU)(0,\tau)
-\nabla_VH_i(\tau,V)\cdot V^{''}(\tau)
-V^{'}(\tau)\cdot \Bigl(\nabla_V(\d_tH_i)(\tau,V)
$$
$$
-(\d_t^2H_i)(\tau,V)-
\nabla_V(\d_tH_i)(\tau,V)
+\nabla_V(\nabla_VH_i(\tau,V)\cdot V^{'}(\tau))\Bigr)
\biggr]
$$
$$
+\Lambda_i^{-1}(0,\tau)\Bigl[-(\d_x\Lambda_i)(0,\tau)(\d_xU_i)(0,\tau)+
(\d_xF_i)(0,\tau,U)
$$
$$
+\nabla_UF_i(0,\tau,U)
\cdot(\d_xU)(0,\tau)\Bigr]\Big|_{\tau=t_i(x,t)}\quad
\mbox{if}\,\,t_i(x,t)>0, 1\le i\le k.
$$
The expressions for $(R_{ixx}^{''}U)(x,t)$ if $t_i(x,t)>0$ and $k+1\le i\le
 n$
are similar. The expressions for
$V_i^{''}$ are derived from the system~(\ref{eq:1}), using
suitable differentiations:
$$
V_i^{''}=\d_t^2U_i=
-\d_t\Lambda_i\d_xU_i-\Lambda_i\d_x\d_tU_{i}+\d_tF_i+\nabla_UF_i\cdot\d_tU
$$
$$
=-\d_t\Lambda_i\d_xU_i+\Lambda_i
\Bigl(\d_x\Lambda_i\d_xU_i+\Lambda_i\d_x^2
 U_i-\d_xF_i-\nabla_UF_i\cdot \d_xU\Bigr)
+\d_tF_i+\nabla_UF_i\cdot \d_tU,
$$
where $\d_x^2U$ satisfies~(\ref{eq:U_xx}) and
the right hand side is considered
restricted to $x=0$ if
$1\le i\le k$ and to $x=l$ if $k+1\le i\le n$.
Take $t_2$ satisfying the condition~(\ref{eq:22_0}) with $m=2$
 and the inequality
$t_2<1/q_2$, where
$$
q_2=(nL_F^{\max}+2E_{\Lambda}(1,0;T))(1+nL_H^{\max}).
$$
This shows that the operator
defined by the right hand side of (\ref{eq:U_xx}) has the
contraction property
with respect to the domain
$\overline\Pi^{t_2}$.
 We therefore have the following global a priori estimate:
\begin{equation}\label{eq:U2g}
\begin{array}{ccccc}
\displaystyle
E_U(2,0;T)\le
P_{1,2}\biggl(\frac{1}{1-q_2t_2},n,L_H^{\max}\biggr)\\[4mm]
\displaystyle
\times
P_{2,2}\biggl(n,\max\limits_{x\in[0,l],1\le i\le n}|A_{i}^{''}(x)|,
\max\limits_{0\le \al_1+\al_2\le 2}E_{\Lambda}(\al_1,\al_2;T),
\max\limits_{0\le\al_1+\al_2\le 1}E_{\Lambda^{-1}}(\al_1,\al_2;T),
\\[5mm]
\displaystyle
\max\limits_{1\le|\beta|+\al_1+\al_2\le 2}
E_{\d_U^{|\beta|}F}(\al_1,\al_2),
\max\limits_{1\le|\beta|+\al_1\le 2}E_{\d_V^{|\beta|}H}(\al_1),
L_F^{\max},L_H^{\max},
\max\limits_{\al_1+\al_2=1}E_U(\al_1,\al_2;T)
\biggr),
\end{array}
\end{equation}
where $\beta=(\beta_1,\dots,\beta_n)$,
$\beta_i\in\N_0$, $|\beta|=\beta_1+\dots+\beta_n$,
 $P_{1,2}$ is a polynomial of degree $3\lceil T/t_2\rceil$
with all coefficients identically equal to 1 and $P_{2,2}$
is a polynomial of the third degree with positive constant
coefficients depending only on $T$.

The estimates for $E_U(1,1;T)$ and $E_U(0,2;T)$ now follow from
system~(\ref{eq:1}) and its suitable differentiations.
The claim is proved.

We further proceed by induction. Assume that the problem~(\ref{eq:integral})
has  a unique $\Con^{m-1}(\overline\Pi^T)$-solution
 for an arbitrary $m\in\N$.

{\it Claim 4. The problem~(\ref{eq:1})--(\ref{eq:3}) has  a unique
$\Con^m(\overline\Pi^T)$-solution.}
Similarly to the proof of  Claims 1--3, we
consider the
problem for $\d_x^mU$ (using suitable differentiations
 and integrations of~(\ref{eq:1})--(\ref{eq:3})).
Taking into account the induction
assumption, we apply the contraction
mapping principle to the operator defined by the right-hand side
of the problem for $\d_x^mU$.  It is not difficult to prove
that the operator has the
contraction property with respect to $\overline\Pi^{t_m}$ for some $t_m>0$
satisfying the condition~(\ref{eq:22_0}) and the inequality
$t_m<1/q_m$, where
\begin{equation}\label{eq:q_m}
q_m=(nL_F^{\max}+mE_{\Lambda}(1,0;T))(1+nL_H^{\max}).
\end{equation}
This implies the existence and the uniqueness of
a $\Con_{x,t}^{m,0}(\overline\Pi^{t_m})$-solution to
the problem~(\ref{eq:1})--(\ref{eq:3}). Iterating this local
existence-uniqueness result in domains
$$
(\Pi^{jt_m}\cap\Pi^T)\setminus\overline{\Pi}^{(j-1)t_m}, \quad
1\le j\le \lceil T/t_m\rceil,
$$
we complete the proof of the claim in
$\lceil T/t_m\rceil$ number of steps. In parallel, we arrive at
 the following global estimate:
\begin{equation}\label{eq:Umg}
\begin{array}{ccccc}
\displaystyle
E_U(m,0;T)\le
P_{1,m}\biggl(\frac{1}{1-q_mt_m},n,L_H^{\max}\biggr)\\[4mm]
\displaystyle
\times
P_{2,m}\biggl(n,\max\limits_{x\in[0,l],1\le i\le n}|A_{i}^{(m)}(x)|,
\max\limits_{0\le \al_1+\al_2\le m-1}
E_{\Lambda^{-1}}(\al_1,\al_2;T),
\\[5mm]
\displaystyle
\max\limits_{0\le \al_1+\al_2\le m}E_{\Lambda}(\al_1,\al_2;T),
\max\limits_{1\le|\beta|+\al_1+\al_2\le m}E_{\d_U^{|\beta|}F}(\al_1,\al_2),
\max\limits_{1\le|\beta|+\al_1\le m}
E_{\d_V^{|\beta|}H}(\al_1),\\[4mm]
\displaystyle
L_F^{\max},L_H^{\max},
\max\limits_{1\le\al_1+\al_2\le m-1}E_U(\al_1,\al_2;T)\biggr).
\end{array}
\end{equation}
Here $P_{1,m}$ is a polynomial of degree $3\lceil T/t_m\rceil$
with all coefficients identically equal to 1. Furthermore, $P_{2,m}$
is a polynomial whose degree depends on $m$ but neither on $T$ nor on $t_m$
and whose coefficients are positive constants
 depending only on $m$ and $T$.
The existence and uniqueness of
a $\Con_{x,t}^{\al_1,\al_2}(\overline\Pi^{T})$-solution
where $\al_1+\al_2=m$, now  follow  from the system~(\ref{eq:1}) and
its suitable differentiations. The respective global a priori estimates
for $E_U(\al_1,\al_2;T)$ one can easily  obtain from
 the inequality~(\ref{eq:Umg}) and the induction assumption.
The claim is proved.

The classical smooth solution to the~problem~(\ref{eq:1})--(\ref{eq:3})
satisfying estimates~(\ref{eq:Umg}) in
 $\overline\Pi^T$  for any $m\in\N_0$
can be constructed
by the sequential approximation method.
We now use this solution to
construct a representative of the Colombeau solution.
According to the assumptions of the theorem, we consider all
the initial data as elements of the corresponding Colombeau algebras.
We choose representatives $\lambda$, $a$, $f$, $h$, $L_f$, and $L_h$
of $\Lambda$, $A$, $F$,  $H$, $L_F$, and $L_H$, respectively,
with the properties required in the theorem.
Let $\phi=\varphi\otimes\varphi\in \A_0(\R^2)$.
Consider a prospective representative  $u=u(\phi,x,t)$ of $U$
which is the classical smooth solution to the
problem~(\ref{eq:1})--(\ref{eq:3})
with the initial data $\lambda(\phi,x,t)$, $a(\varphi,x)$,
$f(\phi,x,t,u(\phi,x,t))$,
$h(\varphi,t,v(\varphi,t))$, $L_f(\phi,x,t)$, $L_h(\varphi,t)$, where
$v(\varphi,t)=(u_1(\phi,0,t),\dots,u_{k}(\phi,0,t)$,
$u_{k+1}(\phi,l,t),\dots,$
$u_{n}(\phi,l,t))$.
For the existence part of the proof, we have to
show that $u\in\EM$, i.e. to obtain
moderate growth estimates
of $u(\phi_{\varepsilon},x,t)$ in terms of the
regularization parameter~$\varepsilon$.

Fix $N\in\N$ to be so large that for all $\varphi\in\A_{N}(\R)$
there exists $\eps_0$ such that for all $\eps<\eps_0$
the following conditions are  true:

a) The moderate estimate (see the definition of $\EM$) holds for
$a(\varphi_{\varepsilon},x)$,
$f(\phi_{\varepsilon},x,t,0)$,
  $h(\varphi_{\varepsilon},t,0)$, $L_f(\phi_{\varepsilon},x,t)$,
and $L_h(\varphi_{\varepsilon},t)$.

b) The invertibility estimate (see Theorem~\ref{thm:invert})
holds for $\lambda(\phi_{\varepsilon},x,t)$.

c) The local-$\gamma$-growth
estimate (see Definition~\ref{defn:ga} ) holds
for $\lambda(\phi_{\varepsilon},x,t)$,
 $L_f(\phi_{\varepsilon},x,t)$, and
$L_h(\varphi_{\varepsilon},t)$.

d) The local-$\gamma_1$-growth
estimate holds for $\partial_x
\lambda(\phi_{\varepsilon},x,t)$.

Fix $\varphi\in\A_{N}(\R)$. Let $p_{1,m}(\varphi)$, $p_{2,m}(\varphi)$,
and $q_m(\varphi)$ denote the value of, respectively, $P_{1,m}$,
$P_{2,m}$, and $q_m$,
 where $U(x,t)$, $\Lambda(x,t)$,  $A(x)$, $F(x,t,U(x,t))$,
 $H(t,V(t))$, $L_F(x,t)$, and $L_H(t)$ are replaced by their representatives $u(\phi,x,t)$,
$\lambda(\phi,x,t)$, $a(\varphi,x)$,
$f(\phi,x,t,$ $u(\phi,x,t))$,
$h(\varphi,t,v(\varphi,t))$, $L_f(\phi,x,t)$, and
$L_h(\varphi,t)$,
respectively.
It suffices to prove the moderate estimates for $p_{1,m}(\varphi)$ and
$p_{2,m}(\varphi)$ for all $m\in\N_0$.
The expression~(\ref{eq:q_m}) and
assumptions imposed on
$\Lambda$, $F$, and $H$
 make it obvious now that
$$
q_m(\varphi_{\eps})\le \gamma^{2N+1}(\varepsilon)+\gamma_1^{2N}(\varepsilon)
$$
for all sufficiently small  $\varepsilon$.
Since $t_m\le\min\{L/
E_{\Lambda}(0,0),1/q_m\}$
and $E_{\Lambda}(0,0)$
$\le \gamma^{N+1}(\varepsilon)$,
we can choose
$t_m=1/[2(\gamma^{2N+1}(\varepsilon)+\gamma_1^{2N}(\varepsilon))]$.
Taking into account~(\ref{eq:ga}), for each $m\in\N_0$ and for all small enough $\varepsilon$ we have
$$
\begin{array}{ccccc}
\displaystyle
\biggl(\frac{nL_H^{\max}}{1-q_mt_m}\biggr)^{3\lceil T/t_m\rceil}\le
\gamma(\varepsilon)^{6(N+1)\lceil T(\gamma^{2N+1}(\varepsilon)+\gamma_1^{2N}
(\varepsilon))\rceil}
\nonumber\\[6mm]
\displaystyle
\le

\Bigl(\gamma(\varepsilon)^{\gamma^{2N+2}(\varepsilon)}\Bigr)
\Bigl(\gamma(\varepsilon)^{\gamma_1^{2N+1}(\varepsilon)}\Bigr)
=O\biggl(
\frac{1}{\varepsilon^2}\biggr)\quad \mbox{as}\,\,\eps\to 0.\nonumber
\end{array}
$$
We conclude that for each $m\in\N_0$ there exists $N\in\N$ such that for all
$\varphi\in\A_{N}(\R)$ we have
\begin{equation}\label{eq:P1m}
\begin{array}{cc}
\displaystyle
P_{1,m}(\varphi_{\eps})
=O\biggl(\frac{1}{\varepsilon^2}\biggr)\quad \mbox{as}\,\,\eps\to 0.
\end{array}
\end{equation}
One can easily see now that for
$\al_1=\al_2=0$
\begin{equation}\label{eq:dxU}
E_{u_{\eps}}(\al_1,\al_2;T)
=O\biggl(\frac{1}{\varepsilon^N}\biggr)\quad \mbox{as}\,\,\eps\to 0
\end{equation}
 for all $\varphi\in\A_{N}(\R)$ with large enough $N\in\N$,
where $u_{\eps}=u(\phi_{\eps},x,t)$.
To prove similar estimates for all derivatives of $U$,
we use induction on $\al=\al_1+\al_2$. Assuming~(\ref{eq:dxU})  to hold for
$\al\le m-1$, let us show that~(\ref{eq:dxU})  is true for $\al=m$ as well.
Indeed, let $\varphi\in\A_{N}(\R)$ with $N$
chosen so large that for all sufficiently small $\eps$ the following
conditions are true:

a) The moderate estimate holds for  $\partial_x^{\al_1}
\partial_t^{\al_2}u_i(\phi_{\varepsilon},x,t)$ for
$0\le \al_1+\al_2\le m-1$ and $1\le i\le n$
(the induction assumption).

b) The moderate estimate holds
for
$\max\limits_{0\le \al_1+\al_2\le m}E_{\lambda_{\eps}}(\al_1,\al_2;T)$,
$\max\limits_{0\le\al_1+\al_2\le m-1}
E_{\lambda_{\eps}^{-1}}(\al_1,\al_2;T)$,
e$a^{(m)}(\varphi_{\varepsilon},x)$,
where $\lambda_{\eps}=\lambda(\phi_{\eps},x,t)$.

c) Given
$U=[(u(\phi_{\eps},x,t))_{\phi\in\A_0(\R^2)}]\in\G(\overline\Pi^T)$ satisfying the estimate~(\ref{eq:dxU}) for
$\al_1=\al_2=0$, the moderate estimate holds for
$\max\limits_{1\le|\beta|+\al_1+\al_2\le m}E_{\d_U^{|\beta|}f_{\eps}}(\al_1,\al_2)$,
$\max\limits_{1\le|\beta|+\al_1\le m}E_{\d_V^{|\beta|}h_{\eps}}(\al_1)$,
where $f_{\eps}=f(\phi_{\eps},x,t,u_{\eps})$,
$h_{\eps}=h(\varphi_{\eps},t,v(\varphi_{\eps},t))$

d) the invertibility
estimate holds for $\lambda(\phi_{\varepsilon},x,t)$.\\
Since $p_{2,m}(\varphi_{\eps})$ is a polynomial whose degree does not depend on $\varepsilon$,
the moderateness of $p_{2,m}(\varphi)$ becomes obvious.
The moderateness of $E_u(m,0;T)$ are done by~(\ref{eq:P1m}).
The moderateness property of $E_u(\al_1,\al_2;T)$ for all other
$\al_1$ and $\al_2$
such that $\al_1+\al_2=m$ is a consequence of the moderateness of
$E_u(m,0;T)$,
the system~(\ref{eq:1}), its suitable differentiations,
and the induction assumption.

Since $T>0$ is arbitrary, the existence part of the proof is complete.

The proof of the uniqueness part  follows the same scheme.
The only difference is that now we
consider the problem with respect to the difference
$U-W$ of two Colombeau solutions $U$ and
$W$. We hence have the problem~(\ref{eq:1})--(\ref{eq:3}) with the right hand sides
$$
\int\limits_0^1\nabla_UF(x,t,\sigma U+(1-\sigma)W)\,d\sigma\cdot(U-W)+M_1,
$$
$$
\int\limits_0^1\nabla_VH(t,\sigma V+(1-\sigma)V_W)\,d\sigma\cdot(V-V_W)+M_3,
$$
and $M_2$ in~(\ref{eq:1}), (\ref{eq:3}), and~(\ref{eq:2}), respectively.
Here $M_i\in\NN$ and $V_W$ is equal to $V$, where $U$ is
replaced by $W$. The analysis is even simpler since,
due to~\cite{gro},  it suffices to check the negligibility
of $U-W$ at order zero. For this purpose we rewrite the
estimate~(\ref{eq:240})  with respect to the function $U-W$ and  use
Assumptions~7--9 and the fact that
$[(p_{2,0}(\varphi))_{\varphi\in\A_0(\R)}]\in\NN$. This finishes the proof.
\end{proof}

\subsection{Non-Lipshitz nonlinearities}

We here extend the above existence-uniqueness result to the case of
non-Lipshitz nonlinearities in~(\ref{eq:1}) and~(\ref{eq:3}).
Set
$$
E_{\nabla F}(U)=\max\Bigl\{|
\nabla_U F_i(x,t,U(x,t))|\,:\,
(x,t)\in\overline{\Pi}^T,1\le i\le n\Bigr\},
$$
$$
E_{\nabla H}(V)=\max\Bigl\{|
\nabla_V H_i(t,V(t))|\,:\,
t\in[0,T],1\le i\le n\Bigr\}.
$$
Simplifying the notation, we drop the dependence of
$E_{\nabla F}(U)$ and $E_{\nabla H}(V)$ on $T$.
Note that we will use these parameters for a fixed $T>0$.

To state the main result of this section, we suppose that at least
one of the following two assumptions holds.
\vskip5mm
{\it Assumption 11.}

a) $H(t,V)$ is smooth in $t,V$ and the mapping $V\mapsto\nabla_VH(t,V)$ is
globally bounded, uniformly over $t$ varying in compact subsets of
$[0,\infty)$;

b) Given $T>0$, there exists $C_F$ such that for all $1\le i\le n$,
$(x,t)\in\overline\Pi^T$, and $y\in\R^n $ we have
$$
|\nabla_U F_i(x,t,y)|\le C_F\log\log D(x,t,y),
$$
where $D(x,t,y)$ with respect to $y$ is a polynomial with
coefficients in $\G(\overline\Pi^T)$.
\vskip5mm
{\it Assumption 12.}

a) Given $T>0$, there exists $C_H$ such that for all $1\le i\le n$,
$t\in[0,T]$, and $z\in\R^n$ we have
$$
|\nabla_V H_i(t,z)|\le C_H(\log\log B(t,z))^{1/4},
$$
where $B(t,z)$  with respect to $z$ is a polynomial with
coefficients in $\G[0,T]$.

b) Assumption~3 is true with $\gamma_1(\eps)=O((\log\log 1/\eps)^{1/4})$;

c) Given $T>0$, there exists $C_F$ such that for all $1\le i\le n$,
$(x,t)\in\overline\Pi^T$, and $y\in\R^n$  we have

$$
|\nabla_U F_i(x,t,y)|\le C_F
(\log\log D(x,t,y))^{1/4},
$$
where $D(x,t,y)$ with respect to $y$  is a polynomial with
coefficients in $\G(\overline\Pi^T)$.

\begin{thm}
Assume that Assumption~11 or~12 is true.
Under Assumptions 1--6 and 10 where the functions $\gamma$ and $\gamma_1$
are specified by~(\ref{eq:ga}), the
problem~(\ref{eq:1})--(\ref{eq:3})
has a unique solution $U\in\G(\overline{\Pi})$.
\end{thm}

\begin{proof}
In the proof we will use a modified notion of $\EM(\overline\Pi)$,
namely, let $u\in\EM(\overline\Pi)$ iff $u\in\E(\overline\Pi)$ and
for every compact set $K\subset\overline\Pi$
there is $N\in\N$ such that for every $\varphi\in \A_{N}(\R^n)$
there exists $\eta>0$ with
$\sup\limits_{x\in K}|u(\varphi_{\varepsilon}\times\varphi_{\varepsilon},x,t)|
\le \gamma^{N}(\varepsilon)$ for
 all
$0<\varepsilon<\eta$.

Fix an arbitrary $T>0$.
From the proof of Theorem~\ref{thm:gen}
it follows that
$[(p_{2,0}(\varphi))_{\varphi\in\A_0(\R)}]$ is a Colombeau generalized number and
hence has the moderateness property.
This means that there
exists $N_1\in\N$ such that for all $\varphi\in\A_{N_1}(\R)$
there is $\eta(\varphi)>0$ with
\begin{equation}\label{eq:mod}
|p_{2,0}(\varphi_{\eps})|\le \eps^{-N_1},\quad 0<\eps<\eta(\varphi).
\end{equation}
Without loss of generality we can assume that $N_1$ is so large that
for all $\varphi\in\A_{N_1}(\R)$ the zero-order moderateness 
property holds for
the coefficients of the polynomial  $D(x,t,y)$ (if Assumption~11 
is fulfilled) or for
the coefficients of the polynomials  $D(x,t,y)$ and $B(t,z)$ (if Assumption~12
is fulfilled).
To simplify notation, we can suppose that, given
$\varphi\in\A_{N_1}(\R)$, the value of $\eta(\varphi)$ in~(\ref{eq:mod})
is so small that the zero-order 
moderate estimates for the
coefficients of $D$ and/or $B$ are true for all $\eps<\eta(\varphi)$.
Note that any $U\in\G(\overline{\Pi})$  has
the following property:  there
exists $N_2\in\N$ such that for all $\varphi\in\A_{N_1+N_2}(\R)$
there is  $\eps_0(\varphi)\le\eta(\varphi)$,
 where the value of $\eta(\varphi)$
is the same as in~(\ref{eq:mod}),
 with
\begin{equation}\label{eq:*}
\sup\limits_{\overline\Pi^T}|u(\varphi_{\eps}\times\varphi_{\eps},x,t)|\le
\eps^{-N_1-N_2}, \quad 0<\eps<\eps_0(\varphi),
\end{equation}
with the constant  $N_1$ being the same as in~(\ref{eq:mod}). Obviously, any increase of $N_2$ and any decrease of
$\eps_0(\varphi)$
will keep this property true. This will allow us to adjust the values of $N_2$ and $\eps_0(\varphi)$ according to our
purposes.

Following the proof of
Theorem~\ref{thm:gen}, for all $\varphi\in\A_{N_1+N_2}(\R)$,  we arrive at the  estimates~(\ref{eq:240})
and~(\ref{eq:Umg})
with $E_{u_{\eps}}(m,0)$, $E_{\nabla f_{\eps}}(u_{\eps})$,
and $E_{\nabla h_{\eps}}(u_{\eps})$
in place of $E_U(m,0)$, $L_F$, and $L_H$, respectively, where $0<\eps<\eta(\varphi)$ and 
 the value of $\eta(\varphi)$
is the same as in~(\ref{eq:mod}). Recall that
$u$, $f$, and $h$ are
representatives of $U$, $F$, and $H$, respectively, and
$u_{\eps}(x,t)=u(\phi_{\eps},x,t)$.
On the account of these estimates,
we will obtain the existence once we
prove the following assertion:
\vskip4mm
$(\iota)$ the constant $N_2\in \N$ can be chosen so that for all
$\varphi\in\A_{N_1+N_2}(\R)$ there exists $\eps_0(\varphi)$ such that
\begin{equation}\label{eq:**}
\Bigl[2n(1+E_{\nabla h_{\eps}}
(u_{\eps}))\Bigr]^
{7T(1+E_{\nabla h_{\eps}}(u_{\eps}))
(E_{\nabla f_{\eps}}(u_{\eps})+
mE_{\Lambda}(1,0;T))}\le \eps^{-N_2},\quad 0<\eps<\eps_0(\varphi),
\end{equation}
whatsoever
$u(\varphi\times\varphi,x,t)\in\E(\overline\Pi^T)$ satisfying
the inequality~(\ref{eq:*}).
\vskip4mm

Let us prove Assertion $(\iota)$ using
Assumption~11. Recall that at this point $N_2$ is a constant whose exact value will be fixed below.
 Fix $\varphi\in\A_{N_1+N_2}(\R)$.
By~(\ref{eq:*}) and Assumption~11, there exists $N_3\in\N$   for which the estimate
$$
E_{\nabla f_{\eps}}(u_{\eps})\le C_F\log\log d(\phi_{\eps},x,t,u_{\eps})\le
C_F\log\log\eps^{-N_3},\quad 0<\eps<\eps_0(\varphi),
$$
is true, where $d$ is a representative of $D$.
Furthermore,
there exist $C_1>1$, $C_2>0$, and $k_1,k_2\in\N$  such that the left hand side of~(\ref{eq:**}) is bounded from
above by
$$
C_1^{C_2(\log\log\eps^{-N_3}+\ga_1(\eps))}\le e^{C_3\log\log\eps^{-N_3}}
\ga(\eps)^{k_1\ga_1(\eps)}
\le e^{\log(\log\eps^{-N_3})^{C_3}}\eps^{-k_2}
$$
$$\le
(\log\eps^{-N_3})^{C_3}\eps^{-k_2}
\le N_3^{\lceil C_3\rceil}\eps^{-\lceil C_3\rceil-k_2},
$$
where $C_3=C_2\log C_1$ and $0<\eps<\eps_0(\varphi)$. 
It is important to note that $C_3$ and $k_2$ can be fixed so that the above estimates hold for all $N_2$ and all
$\varphi$. 
 This makes the values
$N_2=2\lceil C_3\rceil+k_2$ and $\eps_0(\varphi)=\min\{\eta(\varphi),N_3^{-\lceil C_3\rceil}\}$, which we now set up, well defined.
Assertion $(\iota)$  now follows from 
the fact that 
$\varphi$ is an arbitrary function from
$\A_{N_1+N_2}(\R)$.

Let us prove Assertion $(\iota)$ using Assumption~12. Following the
same scheme as above, fix $\varphi\in\A_{N_1+N_2}(\R)$,
where $N_2$ will be specified below.
By~(\ref{eq:*}) and Assumption~12, there exist $N_3,N_4\in\N$
such that the following estimates are true:
$$
E_{\nabla f_{\eps}}(u_{\eps})\le C_F\log\log d(\phi_{\eps},x,t,u_{\eps})\le
C_F\log\log(\eps^{-N_3}),\quad 0<\eps<\eps_0(\varphi),
$$
$$
E_{\nabla h_{\eps}}(u_{\eps})\le C_H\log\log b(\varphi_{\eps},t,u_{\eps})\le
C_H\log\log(\eps^{-N_4}),\quad 0<\eps<\eps_0(\varphi),
$$
where $b$ is a representative of $B$.
 Furthermore,
there exist $C_1>1$ and $C_2>0$ 
 such that
the left hand side of~(\ref{eq:**}) is bounded from above by
$$
\Bigl[C_1\log\log (\eps^{-N_4})\Bigr]^{1/2C_2
(\log\log\eps^{-N_3-N_4})^{1/2}}
$$
$$
\le \exp\Bigl\{{C_2\log(\log(\log
(\eps^{-N_4}))^{C_1})^{1/2}
(\log\log \eps^{-N_3-N_4})^{1/2}}\Bigr\}
$$
$$
\le \exp\Bigl\{{C_2\log(\log\eps^{-N_3-N_4})
^{C_1}}\Bigr\}=
(\log \eps^{-N_3-N_4})^{C_1C_2}
$$
$$
=\Bigl((N_3+N_4)\log\eps^{-1}\Bigr)^{C_1C_2}\le
(N_3+N_4)^{\lceil C_1C_2\rceil}\eps^{-\lceil C_1C_2\rceil},
$$
where $0<\eps<\eps_0(\varphi)$.
Note that $C_1$ and $C_2$ can be fixed so that the above estimates hold for all $N_2$ and all
$\varphi$.
We now set 
$N_2=2\lceil C_1C_2\rceil$ and $\eps_0(\varphi)=\min\{\eta(\varphi),(N_3+N_4)^{-\lceil C_1C_2\rceil}\}$ and this value is well defined.
Assertion $(\iota)$ now follows from
 the fact that $\varphi$ is an arbitrary function in
$\A_{N_1+N_2}(\R)$.

Since $T>0$ is arbitrary, the
existence part of the proof
is complete.

The proof of the uniqueness part  follows the same scheme
(cf.\ also the proof of
Theorem~\ref{thm:gen}).
We apply the estimate~(\ref{eq:240})
to the difference of two generalized solutions
to the problem~(\ref{eq:1})--(\ref{eq:3}). From the
existence part of the proof we see that the first factor in the right-hand side
of~(\ref{eq:240}) has the moderateness propery. Since the second
factor is negligible,
the uniqueness follows.
\end{proof}

\begin{ex}\rm
Consider $n=1$ and
$$
F(x,t,U)=(G_1^2(x,t)+G_2^2(x,t)U^2)^{1/2}\log{\log{(G_3^2(x,t)+
G_4^2(x,t)U^2)}}^{1/2},
$$
where $G_i(x,t)\in\G(\overline\Pi)$.
Then
$$
\d_UF(x,t,U)=\frac{G_2^2U}{(G_1^2+G_2^2U^2)^{1/2}}
\log{\log{(G_3^2+G_4^2U^2)}}^{1/2}
$$
$$
+\frac{G_4^2U(G_1^2+G_2^2U^2)^{1/2}}{\log{(
G_3^2+G_4^2U^2)}^{1/2}
(G_3^2+G_4^2U^2)}.
$$
The function $F(x,t,U)$ is non-Lipshitz and satisfies
Assumption~11(b).
\end{ex}

\begin{rem}
The theorem states that, whatsoever singularity of the initial data
of our problem and whatsoever nonlinearities of $F$ and $H$ allowed
by Assumption~11 (or~12), the
problem~(\ref{eq:1})--(\ref{eq:3})
has a unique solution in the Colombeau algebra $\G(\overline{\Pi})$.
\end{rem}

\end{document}